\documentclass[11pt, draft]{amsart}
\usepackage{amssymb, amstext, amscd, amsmath}
\usepackage[all]{xy}

\newtheorem*{theorem}{The Gauge Invariant Uniqueness Theorem}
\newtheorem*{acknow}{Acknowledgements}

\newcommand{\bT}{{\mathbb{T}}}
\newcommand{\bZ}{{\mathbb{Z}}}

\newcommand{\B}{{\mathcal{B}}}
\newcommand{\F}{{\mathcal{F}}}

\newcommand{\K}{{\mathcal{K}}}
\renewcommand{\L}{{\mathcal{L}}}
\renewcommand{\O}{{\mathcal{O}}}
\newcommand{\T}{{\mathcal{T}}}

\def\be{\beta}
\def\ga{\gamma}

\newcommand{\foral}{\text{ for all }}
\newcommand{\forr}{\text{ for }\ }
\newcommand{\qand}{\quad\text{and}\quad}

\newcommand{\ad}{\operatorname{ad}}
\newcommand{\id}{{\operatorname{id}}}
\newcommand\Span{\mathop{\rm span}}

\newcommand{\ca}{\mathrm{C}^*}
\newcommand{\ol}{\overline}

\newcommand{\sca}[1]{\left\langle#1\right\rangle} 
\newcommand{\nor}[1]{\left\Vert #1\right\Vert}    

\begin{document}

\title[A note on C*-correspondences]{A Note on the Gauge Invariant Uniqueness Theorem for C*-correspondences}

\author[E.T.A. Kakariadis]{Evgenios T.A. Kakariadis}
\address{Department of Mathematics\\Ben-Gurion University of the Negev\\Be'er Sheva 84105\\Israel}
\email{kakariad@math.bgu.ac.il}

\thanks{2010 {\it  Mathematics Subject Classification.}
46L08, 47L55}
\thanks{{\it Key words and phrases:} C*-correspondences, Toeplitz-Cuntz-Pimsner algebras.}

\maketitle

\begin{abstract}
We present a short proof of the gauge invariant uniqueness theorem for relative Cuntz-Pimsner algebras of C*-correspondences.
\end{abstract}

\section{Introduction}

Soon after its initiation by Pimsner \cite{Pim97}, the theory of C*-correspondences captured the interest of the research community. The motivating feature is their flexible language that encodes a broad variety of examples in operator algebras (both selfafdjoint and nonselfadjoint). Nowadays they give a central construction in the general theory of C*-algebras \cite{BroOza08}, following the general framework provided by Katsura \cite{Kat04-1, Kat04-2, Kat07}. In this note we present a short proof of the gauge invariant uniqueness theorem for relative Cuntz-Pimsner algebras by using a less sharp analysis of the cores than that of \cite{Kat04-2}.

The terminology of C*-correspondences has been under a number of considerable changes in the last years. In this paper we follow \cite{Kat04-2}. A \emph{C*-correspondence $X$ over a C*-algebra $A$} is a right Hilbert $A$-module along with a $*$-homomorphism $\phi_X \colon A \to \L(X)$ on the adjointable operators $\L(X)$. We say that a pair $(\pi,t)$ defines a representation of $X$ if $\pi\colon A \to \B(H)$ is a $*$-representation and $t \colon X \rightarrow \B(H)$ is a linear map such that $\pi(a)t(\xi)\pi(b) = t(\phi_X(a)\xi b)$ and $t(\xi)^*t(\eta) = \pi(\sca{\xi,\eta}_X)$ for all $a,b \in A$ and $\xi, \eta \in X$. The C*-property implies that $t$ is isometric when $\pi$ is injective. Kajiwara, Pinzari and Watatani \cite{KPW98} show that an (injective) pair $(\pi,t)$ induces an (injective) $*$-representation $\psi_t \colon \K(X) \to \B(H)$ such that $\psi_t(\theta_{\xi,\eta})= t(\xi)t(\eta)^*$. As we are about to see, the ideal
\[
I_{(\pi,t)}': = \{ a\in A \mid \pi(a) \in \psi_t(\K(X))\},
\]
plays a significant role (see also Katsura's work \cite{Kat04-2, Kat07}). We say that $(\pi,t)$ \emph{admits a gauge action $\{\ga_z\}_{z \in \bT}$}, if $\{\ga_z\}_{z \in \bT}$ is a point-norm continuous family of $*$-endomorphisms with
\begin{align*}
\ga_z(\pi(a)) = \pi(a), \forr a \in A, \quad \ga_z(t(\xi)) = z t(\xi), \forr \xi \in X.
\end{align*}

Let $J$ be an ideal of $A$ contained in $\phi_X^{-1}(\K(X))$. We say that a pair $(\pi,t)$ is \emph{$J$-covariant} if $\pi(a) = \psi_t(\phi_X(a))$ for all $a\in J$. Then $\O(J,X)$ is the universal C*-algebra generated by (the copies of) $A$ and $X$ relative to $J$-covariant pairs $(\pi,t)$. When $J = (0)$ we denote $\O(J,X)$ by $\T_X$ which is called the \emph{Toeplitz-Pimsner algebra}. When $J$ is Katsura's ideal \cite{Kat04-2}
\[
J_X:= \ker\phi_X^\perp \cap \phi_X^{-1}(\K(X))
\]
we write $\O_X \equiv \O(J_X,X)$ which is called the \emph{Cuntz-Pimsner algebra}. In particular $\T_X$ is the universal C*-algebra relative to the representations of $X$, hence $\O(J,X)$ is the quotient of $\T_X$ by the ideal generated by the elements $\pi(a) - \psi_t(\phi_X(a))$ for all $a\in J$.  It should be noted that Pimsner \cite{Pim97} considers C*-algebras generated simply by a copy of $X$. However the C*-correspondences therein are assumed to be (injective and) \emph{full}, i.e., $\sca{X,X}$ is dense in $A$ \cite[Remark 1.2 (3)]{Pim97}. That is the reason why the C*-algebras in \cite{Pim97} that are generated simply by $X$ manage to reconstruct a copy of $A$.

There is a strong connection between $\T_X$ and $\O_X$ attained by Muhly and Solel \cite{MuhSol98} and Fowler, Muhly and Raeburn \cite{FMR03} under certain assumptions, and settled in full generality by Katsoulis and Kribs \cite{KatKri06}. The \emph{tensor algebra} $\T_X^+$ in the sense of Muhly and Solel \cite{MuhSol98} is the non-involutive closed subalgebra of $\T_X$ generated by $A$ and $X$. Then $\O_X$ is the minimal C*-cover of $\T_X^+$, i.e., the \emph{C*-envelope} of $\T_X^+$ in the sense of Arveson \cite{Arv69}.

A key role in the theory of C*-correspondences is played by the gauge invariant uniqueness theorems. This type of result was initiated by an Huef and Raeburn for Cuntz-Krieger algebras \cite[Theorem 2.3]{HueRae97} and various generalizations were given by Doplicher, Pinzari and Zuccante \cite[Theorem 3.3]{DPZ98}, Fowler, Muhly and Raeburn \cite[Theorem 4.1]{FMR03}, and Fowler and Raeburn \cite[Theorem 2.1]{FowRae99}. In all these cases at least injectivity of $\phi_X$ is assumed, therefore they were not enough to treat general constructions such as C*-algebras of graphs with sources. Gauge invariant uniqueness theorems for $\T_X$ and $\O_X$ were given in full generality by Katsura \cite[Theorem 6.2, Theorem 6.4]{Kat04-2} by using a sharp analysis of the ideal structure of the cores and a conceptual argument concerning short exact sequences. An alternative proof for $\O_X$ was given afterwards by Muhly and Tomforde \cite{MuhTom04} by using a tail adding technique. An extended gauge invariant uniqueness theorem in the much broader class of C*-algebras associated to the pre-C*-correspondences is given by Kwa\'{s}niewski \cite{Kwa11}.

\begin{theorem}
Let $X$ be a C*-correspon\-dence over $A$ and let $J$ be an ideal of $A$ contained in $J_X$. Then a pair $(\pi,t)$ defines a faithful representation of the $J$-relative Cuntz-Pimsner algebra $\O(J,X)$ if and only if $(\pi,t)$ admits a gauge action, $\pi$ is injective and $I_{(\pi,t)}'=J$.
\end{theorem}

As an immediate consequence of the gauge invariant uniqueness theorem we obtain that if $(\pi,t)$ admits a gauge action and $\pi$ is injective then $\ca(\pi,t) \simeq \O(I_{(\pi,t)}',X)$. We remark that for Toeplitz-Pimsner algebras the condition $I_{(\pi,t)}'=(0)$ implies injectivity of $\pi$, whereas for the Cuntz-Pimsner algebra $\O_X$ the condition $I_{(\pi,t)}'=J_X$ is redundant. Indeed when $\pi$ is injective then Katsura remarks that $I_{(\pi,t)}' \subseteq J_X$ \cite[Proposition 3.3]{Kat04-2}. Finally, when $J \subseteq \phi^{-1}_X(\K(X))$ then $A$ (and consequently $X$) embeds isometrically in $\O(J,X)$, if and only if $J \subseteq J_X$, if and only if $\phi_X|_J$ is injective \cite[Lemma 2.7]{KakPet13}. Therefore Katsura's ideal $J_X$ is the maximal ideal for obtaining the minimal C*-algebra that contains an isometric copy of $X$.


\section{The Proof}

\subsection*{Preliminaries}

We follow notation and terminology of \cite{Kat04-2}. We write $X^{\otimes 0} = A$ and $X^{\otimes n+1}= X^{\otimes n} \otimes_A X$ for $n \geq 0$, i.e., \emph{the $A$-stabilized tensor product}. Every $X^{\otimes n}$ becomes a C*-correspondence over $A$ by $\phi_{X^{\otimes n}} = \phi_X \otimes \id_{n-1}$ when $n \geq 1$, and $\phi_{X^0}$ is the multiplication action on $A$. We write $\phi_n \equiv \phi_{X^{\otimes n}}$ for simplicity.

Given a representation $(\pi,t)$ of a C*-correspondence $X$ that acts on a Hilbert space $H$, let $\ca(\pi,t)$ be the C*-subalgebra of $\B(H)$ generated by $\pi(A)$ and $t(X)$. We denote by $(\pi,t^n)$ the induced pair on $X^{\otimes n}$ such that $t^n(\xi_1\otimes \dots \otimes \xi_n) = t(\xi_1) \dots t(\xi_n)$ for all $\xi_1, \dots, \xi_n \in X$. Then
\[
\ca(\pi,t)=\ol{\Span}\{ t^n(\ol{\xi})t^m(\ol{\eta})^* \mid \ol{\xi} \in X^{\otimes n}, \ol{\eta} \in X^{\otimes m}, n,m \in \bZ_+ \}.
\]

When $\pi$ is injective then the equation $\pi(a) = \psi_t(k)$ implies that $\phi_X(a)=k$ and $a\in J_X$ \cite[Proposition 3.3]{Kat04-2}. In short, the C*-identity implies that $t$ is injective when $\pi$ is injective, therefore
\begin{align*}
\nor{\phi_X(a)\xi - k\xi}_X
& = \nor{t(\phi_X(a)\xi) - t(k\xi)}_{\B(H)} \\
& = \nor{(\pi(a)-\psi_t(k))t(\xi)}_{\B(H)} =0,
\end{align*}
for all $\xi \in X$, which shows that $\phi_X(a) =k$. Moreover, since $\pi(b)t(\xi)t(\eta)^* = t(\phi_X(b)\xi)t(\eta)^*$, then for $b \in \ker\phi_X$ we obtain
\[
\pi(ba) = \pi(b) \psi_t(k) = \psi_t(\phi_X(b)k)= 0,
\]
which shows that $a \in \ker\phi_X^\perp$ due to the injectivity of $\pi$.

Let the \emph{cores} of $\ca(\pi,t)$ be the C*-subalgebras
\[
B_{[l,m]} = \Span \{\psi_{t^n}(k_n) \mid k_n \in \K(X^{\otimes n}), l \leq n \leq m\}.
\]
To see that the $*$-subalgebras $B_{[l,m]}$ are indeed closed, first note that $B_{[l,l]}$ is closed since $\psi_{t^l}$ has closed range. Moreover $\psi_{t^{l+1}}(\K(X^{\otimes l+1}))$ is a closed ideal of $B_{[l,l+1]}$ hence
\[
B_{[l,l+1]} = B_{[l,l]} + B_{[l+1,l+1]}= q^{-1} \circ q(B_{[l,l]})
\]
is closed, where $q \colon  \ol{B_{[l,l+1]}} \to \ol{B_{[l,l+1]}} / B_{[l+1,l+1]}$ is the usual quotient $*$-epimorphism. Inductively we get that $B_{[l,m]}$ is a C*-subalgebra of $\ca(\pi,t)$.

If $(e_i)$ is an approximate identity of $\K(X)$ then $(\psi_t(e_i))$ is an approximate identity of $B_{[n,n]}$ for all $n \geq 1$, since $\psi_{t^n}(\K(X^{\otimes n}))$ is the closure of the linear span of $t(\xi_1)\dots t(\xi_n)t(\eta_n)^*\dots t(\eta_1)^*$ for $\xi_i, \eta_i \in X$. Consequently $(\psi_t(e_i))$ is an approximate identity for $B_{[1,m]}$ for all $m \geq 1$. Furthermore
\[
t(X)^* \cdot B_{[1,m]} \cdot t(X) \subseteq B_{[0,m-1]}, \foral m \in \bZ_+.
\]

Furthermore a gauge action $\{\ga_z\}_{z \in \bT}$ on a pair $(\pi,t)$ defines the conditional expectation
\[
E(f) = \int_z \ga_z(f) dz, \foral f \in \ca(\pi,t), \foral f \in \ca(\pi,t).
\]
Then the fixed point algebra $E(\ca(\pi,t)) = \ca(\pi,t)^\ga$ is the inductive limit of the C*-subalgebras $B_{[0,n]}$.

\subsection*{The Fock representation}

For $\xi \in X$ let $\tau_n(\xi) \in \L(X^{\otimes n}, X^{\otimes n+1})$ such that $\tau_n(\xi)(\eta_1 \otimes \dots \otimes \eta_n) = \xi \otimes \eta_1 \otimes \dots \otimes \eta_n$, for $\eta_1, \dots, \eta_n \in X$. Let the \emph{Fock space} $\F(X) = \oplus_{n \in \bZ_+} X^{\otimes n}$ be the direct sum Hilbert $A$-module of $X^{\otimes n}$. The \emph{full Fock representation} is then defined by $(\pi,t)$ with
\[
\pi(a) = \sum_{n \geq 0} \phi_n(a), \forr a \in A, \qand t(\xi) = \sum_{n \geq 0} \tau_n(\xi), \forr \xi \in X.
\]
Then $(\pi,t)$ defines a representation of $X$. A useful fact is that $\psi_t \colon \K(X) \rightarrow \L(\F(X))$ takes up the form
\[
\psi_t(k) = \sum_{n \geq 1} k \otimes \id_{n-1}, \foral k \in \K(X).
\]
Indeed it suffices to note that if $k = \theta_{\xi,\eta}$, then $\tau_n(\xi) \tau_n(\eta)^* = \theta_{\xi,\eta} \otimes \id_{n-1} = k \otimes \id_{n-1}$. Therefore for $a \in A$ such that $\phi_X(a)=k \in \K(X)$ we obtain
\[
\pi(a) - \psi_t(\phi(a)) = \phi_0(a) + \sum_{n \geq 1} \left(\phi_n(a) - k \otimes \id_{n-1}\right) = \phi_0(a).
\]
Moreover $(\pi, t)$ admits a gauge action $\{\ga_z\}_{z \in \bT}$ by letting $\ga_z = \ad_{u_z}$ where
\[
u_z(\xi_1 \otimes \dots \otimes \xi_n) := z^n \xi_1 \otimes \dots \otimes \xi_n,
\]
defines an adjointable unitary operator in $\L(\F(X))$.

Let $J$ be an ideal of $A$ contained in $\phi_X^{-1}(\K(X))$. Then
\[
\K(XJ) = \ol{\Span}\{\theta_{\xi a, \eta} \mid \xi, \eta \in X, a \in J\},
\]
is a closed ideal in $\K(X)$. A crucial remark is that if $\phi_X(a) \in \K(X)$ then $\phi_X(a) \in \K(XJ)$ if and only if $\sca{\eta,\phi_X(a)\xi}_X \in J$, which follows by \cite[Lemma 2.6]{FMR03} or \cite[Lemma 1.6]{Kat07}. In short, if $\sca{k\xi,k\xi}_X \in J$ then there is a $\xi' \in X$ and a positive $a \in J$ such that $k\xi = \xi' a$ by \cite[Lemma 4.4]{Lan95}. Then
\[
k \theta_{\xi,\eta} = \theta_{\xi'a, \eta} = \theta_{\xi'\sqrt{a},\eta \sqrt{a}} \in \K(XJ).
\]
Since $k \in \K(X)$ then $k = \lim_i k e_i$ for some approximate identity $(e_i)$ in $\K(X)$ and the above remark shows that the convergent net $(ke_i)$ is in the closed ideal $\K(XJ)$, thus $k \in \K(XJ)$. Conversely if $k$ is the norm limit of some $k_i = \sum_{m=0}^{N_i} \theta_{\xi_m a_m, \eta_m } \in\K(XJ)$ then
\[
\sca{\eta, \theta_{\xi_m a_m, \eta_m } \xi}_X = \sca{\eta,\xi_m}_X a_m \sca{\eta_m, \xi}_X \in J,
\]
which implies that $\sca{\eta,\phi_X(a)\xi}_X \in J$.

Let the quotient $*$-epimorphism $q_J \colon \L(\F(X)) \to \L(\F(X))/\K(\F(X)J)$ where
\[
\K(\F(X)J) = \ol{\Span}\{\theta_{\ol{\xi} a, \ol{\eta}} \mid \ol{\xi}, \ol{\eta} \in \F(X), a \in J\}.
\]
Then $(q_J \circ \pi,q_J \circ t)$ is a well defined $J$-covariant representation, since
\[
\pi(a) - \psi_t(\phi(a)) = \phi_0(a) \, \in \, J \subseteq \K(\F(X)J),
\]
for all $a \in J$. In particular note that $\psi_{q_J \circ t} = q_J \circ \psi_t$ on $\K(X)$. Furthermore $\K(\F(X)J) \subseteq \ca(\pi,t)$ since
\[
\theta_{\ol{\xi} a, \ol{\eta}} = t^n(\ol{\xi}) \phi_0(a) t^m(\ol{\eta})^* = t^n(\ol{\xi}) \left(\pi(a) - \psi_t(\phi(a))\right) t^m(\ol{\eta})^*,
\]
for all $\ol{\xi} \in X^{\otimes n}, \ol{\eta} \in X^{\otimes m}$ and $a \in J$. Note that $\K(\F(X)J)$ is $\ga_z$-invariant for all $z\in \bT$, therefore $(q_J\circ \pi, q_J \circ t)$ inherits the gauge action $\{q_J \circ \ga_z\}_{z \in \bT}$. We will say that $(q_J\circ \pi, q_J \circ t)$ induces the \emph{$J$-relative Fock representation}.

When $J \subseteq J_X:= \ker\phi_X^\perp \cap \phi_X^{-1}(\K(X))$ then $(q_J \circ \pi, q_J \circ t)$ is isometric. This follows as in \cite[Proposition 4.9]{Kat04-2}. In short, the $*$-homomorphism
\[
\K(X^{\otimes n} J) \, \ni \, k \mapsto k \otimes \id_{n} \, \in \, \L(X^{\otimes n+1})
\]
is injective since the restriction of $\phi_X$ on $J$ is injective, and if
\[
0 = \sca{k \ol{\xi} \otimes \eta_1, k \ol{\xi} \otimes \eta_2}_{X^{\otimes n+1}} = \sca{\eta_1, \phi_X(\sca{k \ol{\xi}, k \ol{\xi}}_{X^{\otimes n}}) \eta_2}_X,
\]
for all $\eta_1, \eta_2 \in X$ then the element $\sca{k\ol{\xi}, k\ol{\xi}}_{X^{\otimes n}}$ of $J$ is also in $\ker\phi_X$. If $\pi(a) \in \K(\F(X)J)$ then $\phi_n(a) \in \K(X^{\otimes n} J)$ for all $n$, hence
\[
\nor{a}_A = \lim_n \nor{\phi_n(a)} = \lim_n \nor{P_n \pi(a) P_n} =0,
\]
since $\lim_n P_n k P_n =0$ for all $k \in \K(\F(X))$, where we write $P_n$ for the projection of $\F(X)$ onto the direct summand $X^{\otimes n}$.

One last property of $(q_J \circ \pi, q_J \circ t)$ is that $I_{(q_J \circ \pi, q_J \circ t)}' = J$. Indeed for $a \in I_{(q_J \circ \pi, q_J \circ t)}'$ let $k \in \K(X)$ such that $q_J\circ \pi(a) = \psi_{q_J \circ t}(k)$. Then injectivity of $q_J \circ \pi$ implies that $\phi_X(a) = k$ and
\[
0 = q_J \circ \pi(a) - \psi_{q_J \circ t}(\phi_X(a)) = q_J(\pi(a) - \psi_t(\phi_X(a))) = q_J(\phi_0(a)),
\]
thus $a \in J$. Conversely $J \subseteq I_{(\pi,t)}'$ by the $J$-covariance of $(\pi,t)$.

\subsection*{The proof}

Fix an ideal $J \subseteq J_X$. We will denote by $(\pi_u,t_u)$ the universal representation of $\O(J,X)$. By the existence of the $J$-relative Fock representation we obtain that $\pi_u$ is isometric. Furthermore, the universal property of $\O(J,X)$ provides the existence of a gauge action $\{\be_z\}_{z \in \bT}$ for $(\pi_u,t_u)$. The cores in $\O(J,X)$ will be denoted by $\B_{[l,m]}$.

Suppose that $(\pi,t)$ is a representation of $X$ that admits a gauge action $\{\ga_z\}_{z \in \bT}$, that $\pi$ is injective and that $I_{(\pi,t)}'=J$. Let $\Phi \colon \O(J,X) \to \ca(\pi,t)$ be the canonical $*$-epimorphism. Then $\Phi \circ \be_z = \ga_z \circ \Phi$, and by a usual C*-argument it suffices to show that the restriction of $\Phi$ to the fixed point algebra $\O(J,X)^\be$ is faithful. Since the fixed point algebra is the inductive limit of the cores $\B_{[0,N]}$ it suffices to show that the kernel of $\Phi$ intersects trivially all $\B_{[0,N]}$.

For $N=0$ we have that $\B_{[0,0]} = \pi_u(A)$ and by assumption $\Phi \circ \pi_u|_{A} = \pi$ is injective. For the inductive step let $N \geq 1$ be the least non-negative integer such that $\ker\Phi \cap \B_{[0,N]} \neq (0)$. Therefore $\ker\Phi \cap \B_{[0,N-1]} =(0)$ and let $f = \pi_u(a) + \sum_{n=1}^N \psi_{t_u^n}(k_n) \neq 0$ for $a \in A$ and $k_n \in \K(X^{\otimes n})$ such that $\Phi(f)=0$. Then $\pi(a) = - \sum_{n=1}^N \psi_{t^n}(k_n)$. Let $(e_i)$ be an approximate identity in $\K(X)$ and compute
\begin{align*}
\lim_i \psi_t(\phi_X(a)e_i)
& = \lim_i \pi(a) \psi_t(e_i) \\
& = - \sum_{n=1}^N \lim_i \psi_{t^n}(k_n)\psi_t(e_i)
 = - \sum_{n=1}^N \psi_{t^n}(k_n),
\end{align*}
where we have used that $(\psi_t(e_i))$ acts as an approximate identity on all $\psi_{t^n}(\K(X^{\otimes n}))$. Therefore the net $\left(\psi_t(\phi_X(a)e_i) \right)$ converges and so it is Cauchy in $\B(H)$. Thus so is the net $(\phi_X(a)e_i)$ in $\K(X)$ since $\pi$, and consequently $\psi_t$, is injective. Hence $(\phi_X(a)e_i)$ converges to some compact operator, say $k \in \K(X)$. Therefore
\[
\pi(a) = - \sum_{n=1}^N \psi_{t^n}(k_n) = \lim_i \psi_t(\phi_X(a)e_i) = \psi_t(k).
\]
Since $I_{(\pi,t)}' = J$ we obtain that $a\in J$. Consequently $\pi_u(a) = \psi_{t_u}(\phi_X(a)) \in \B_{[1,N]}$, which implies that $f \in \B_{[1,N]}$. However in this case
\[
t_u(\eta)^* \cdot f \cdot t_u( \xi ) = \sum_{n=1}^{N} t_u(\eta)^* \psi_{t_u^n}(k_n) t_u(\xi) \, \in \, \ker\Phi \cap \B_{[0,N-1]},
\]
for all $\xi, \eta \in X$. By the choice of $N$ we then obtain that $t_u(\eta)^* \cdot f \cdot t_u(\xi)=0$. In particular $\psi_{t_u}(\theta_{\eta_1,\eta_2}) \cdot f \cdot \psi_{t_u}(\theta_{\xi_1,\xi_2}) =0$, for all $\eta_1, \eta_2, \xi_1, \xi_2 \in X$, hence
\[
\psi_{t_u}(\K(X)) \cdot f \cdot \psi_{t_u}(\K(X)) = (0).
\]
But $\psi_{t_u}(\K(X))$ contains an approximate identity for $\B_{[1,N]}$, hence $f=0$ which is a contradiction.

Conversely, suppose that $(\pi,t)$ defines a faithful representation $\Phi$ of the relative $\O(J,X)$ such that $\Phi(\pi_u(a)):= \pi(a)$ and $\Phi(t_u(\xi))=t(\xi)$. Then $(\pi,t)$ admits a gauge action $\{\ga_z\}_{z \in \bT}$ by $\ga_z:=\Phi \circ \be_z \circ \Phi^{-1}$, and $\pi = \Phi \circ \pi_u$ is injective. Finally we obtain that $I_{(\pi,t)}' = J$. Indeed, since $\pi_u(a) = \psi_{t_u}(\phi_X(a)) = \psi_{t_u}(k)$ for $a\in J$, then
\[
\pi(a) = \Phi(\pi_u(a)) = \Phi(\psi_{t_u}(\phi_X(a))) = \Phi(\psi_{t_u}(k)) = \psi_t(k) = \psi_t(\phi_X(a)).
\]
That is $(\pi,t)$ is $J$-covariant, hence $J \subseteq I_{(\pi,t)}'$. Conversely if $\pi(a) = \psi_t(k)$ for some $k \in \K(X)$ then injectivity of $\pi$ implies that $\phi_X(a)=k$ and
\[
\pi_u(a) = \Phi^{-1}(\pi(a)) = \Phi^{-1}(\psi_t(k)) = \psi_{t_u}(\phi_X(a)).
\]
In particular this holds for the injective $J$-relative Fock representation, thus $a \in J$, and the proof is complete.

\begin{acknow}
\textup{The author was partially supported by the Kreitman Foundation Post-doctoral Fellow Scholarship, and  by a postdoctoral fellowship funded by the Skirball Foundation via the Center for Advanced Studies in Mathematics at Ben-Gurion University of the Negev.}
\end{acknow}


\end{document}